\documentclass[a4paper,11pt]{article}

\usepackage{amsfonts,amssymb,amsmath,amsthm,latexsym,epsf,epsfig,amscd,bbm,stmaryrd}

\usepackage[english]{babel}
\makeatletter \@addtoreset{equation}{section} \makeatother
\makeatletter \@addtoreset{enunciato}{section} \makeatother
\newcounter{enunciato}[section]

\newtheorem{ittheorem}{Theorem}
\newtheorem{itlemma}{Lemma}
\newtheorem{itproposition}{Proposition}
\newtheorem{itdefinition}{Definition}
\newtheorem{itremark}{Remark}
\newtheorem{itclaim}{Claim}
\newtheorem{itfact}{Fact}
\newtheorem{itconjecture}{Conjecture}
\newtheorem{itcorollary}{Corollary}

\newenvironment{theorem}{\addtocounter{enunciato}{1}
\begin{ittheorem}}{\end{ittheorem}}
\newenvironment{lemma}{\addtocounter{enunciato}{1}
\begin{itlemma}}{\end{itlemma}}
\newenvironment{proposition}{\addtocounter{enunciato}{1}
\begin{itproposition}}{\end{itproposition}}
\newenvironment{definition}{\addtocounter{enunciato}{1}
\begin{itdefinition}}{\end{itdefinition}}
\newenvironment{remark}{\addtocounter{enunciato}{1}
\begin{itremark}}{\end{itremark}}

\newenvironment{conjecture}{\addtocounter{enunciato}{1}
\begin{itconjecture}}{\end{itconjecture}}
\newenvironment{corollary}{\addtocounter{enunciato}{1}
\begin{itcorollary}}{\end{itcorollary}}
\newcommand{\be}[1]{\begin{equation}\label{#1}}
\newcommand{\ee}{\end{equation}}
\newcommand{\bl}[1]{\begin{lemma}\label{#1}}
\newcommand{\el}{\end{lemma}}
\newcommand{\br}[1]{\begin{remark}\label{#1}}
\newcommand{\er}{\end{remark}}
\newcommand{\bt}[1]{\begin{theorem}\label{#1}}
\newcommand{\et}{\end{theorem}}
\newcommand{\bd}[1]{\begin{definition}\label{#1}}
\newcommand{\ed}{\end{definition}}
\newcommand{\bp}[1]{\begin{proposition}\label{#1}}
\newcommand{\ep}{\end{proposition}}
\newcommand{\bc}[1]{\begin{corollary}\label{#1}}
\newcommand{\ec}{\end{corollary}}
\newcommand{\bcj}[1]{\begin{conjecture}\label{#1}}
\newcommand{\ecj}{\end{conjecture}}

\def \Z {{\mathbb Z}}
\def \R {{\mathbb R}}

\def \N {{\mathbb N}}

\def \ba {\begin{array}}
\def \ea {\end{array}}

\def \P {{\mathbb P}}
\def \E {{\mathbb E}}
\def \cP {{\mathcal P}}
\def \cE {{\mathcal E}}
\def \cG {{\mathcal G}}

\def \e {{\bar{e}}}
\def \eps {\epsilon}

\def \erom  {\mathrm{e}}
\def \di {\mathrm{d}}
\def \trn {{\interleave}}

\def \SRW {{\hbox{\tiny\rm SRW}}}
\def \RW {{\hbox{\tiny\rm RW}}}

\def\1{\mathbbm{1}}

\begin{document}

\title{Large deviation principle for one-dimensional\\
random walk in dynamic random environment:\\
attractive spin-flips and simple symmetric exclusion}

\author{\renewcommand{\thefootnote}{\arabic{footnote}}
L.\ Avena \footnotemark[1]
\\
\renewcommand{\thefootnote}{\arabic{footnote}}
F.\ den Hollander \footnotemark[1]\,\,\,\,\footnotemark[2]
\\
\renewcommand{\thefootnote}{\arabic{footnote}}
F.\ Redig \footnotemark[1]}

\footnotetext[1]{
Mathematical Institute, Leiden University, 
P.O.\ Box 9512, 2300 RA Leiden, The Netherlands} 

\footnotetext[2]{
EURANDOM, P.O.\ Box 513, 5600 MB Eindhoven, The Netherlands}

\maketitle

\begin{abstract}

Consider a one-dimensional shift-invariant attractive spin-flip system 
in equilibrium, constituting a dynamic random environment, together with 
a nearest-neighbor random walk that on occupied sites has a local drift 
to the right but on vacant sites has a local drift to the left. In \cite{AvdHoRe09} 
we proved a law of large numbers for dynamic random environments satisfying
a space-time mixing property called cone-mixing. If an attractive spin-flip 
system has a finite average coupling time at the origin for two copies starting 
from the all-occupied and the all-vacant configuration, respectively, then it 
is cone-mixing.

In the present paper we prove a large deviation principle for the empirical 
speed of the random walk, both quenched and annealed, and exhibit some properties 
of the associated rate functions. Under an exponential space-time mixing condition 
for the spin-flip system, which is stronger than cone-mixing, the two rate 
functions have a unique zero, i.e., the slow-down phenomenon known to be possible 
in a static random environment does not survive in a fast mixing dynamic random
environment. In contrast, we show that for the simple symmetric exclusion dynamics, 
which is not cone-mixing (and which is not a spin-flip system either), slow-down 
does occur.

\vspace{0.5cm}\noindent

\vspace{0.5cm}\noindent
$\ast$ Invited paper to appear in the 15-th anniversary celebration issue of 
\emph{Markov Processes and Related Fields}. 
\end{abstract}

\newpage


\section{Introduction and main results}
\label{S1}


\subsection{Random walk in dynamic random environment: attractive spin-flips}
\label{S1.1}

Let 
\be{asfsdef} 
\xi = (\xi_t)_{t \geq 0} \quad \mbox{ with } \quad \xi_t=\{\xi_t(x)\colon\,x\in\Z\} 
\ee
denote a one-dimensional \emph{spin-flip system}, i.e., a Markov process on 
state space $\Omega=\{0,1\}^\Z$ with generator $L$ given by 
\be{Generator} 
(Lf)(\eta) = \sum_{x\in\Z} c(x,\eta)[f(\eta^x)-f(\eta)], \qquad \eta\in\Omega,
\ee 
where $f$ is any cylinder function on $\Omega$, $c(x,\eta)$ is the local 
rate to flip the spin at site $x$ in the configuration $\eta$, and $\eta^x$ is 
the configuration obtained from $\eta$ by flipping the spin at site $x$.
We think of $\xi_t(x)=1$ ($\xi_t(x)=0$) as meaning that site $x$ is occupied 
(vacant) at time $t$. We assume that $\xi$ is \emph{shift-invariant}, i.e.,
for all $x\in\Z$ and $\eta\in\Omega$, 
\be{shift}
c(x,\eta)= c(x+y,\tau_y\eta), \qquad y\in\Z,
\ee
where $(\tau_y\eta)(z)=\eta(z-y)$, $z\in\Z$, and also that $\xi$ is \emph{attractive}, 
i.e., if $\eta\leq\zeta$, then, for all $x\in\Z$,
\be{attractiveness}
\begin{aligned}
&c(x,\eta)\leq c(x,\zeta) \quad \text{ if } \quad \eta(x)=\zeta(x)=0,\\
&c(x,\eta)\geq c(x,\zeta) \quad \text{ if } \quad \eta(x)=\zeta(x)=1.
\end{aligned}
\ee
For more on shift-invariant attractive spin-flip systems, see \cite{Li85}, Chapter III. 
Examples are the (ferromagnetic) Stochastic Ising Model, the Voter Model, the Majority 
Vote Process and the Contact Process. 

We assume that
\be{erg}
\xi \mbox{ has an equilibrium } \mu \mbox{ that is shift-invariant and shift-ergodic}.
\ee 
For $\eta\in\Omega$, we write $P^\eta$ to denote the law of $\xi$ starting from 
$\xi(0)=\eta$, which is a probability measure on path space $D_\Omega[0,\infty)$, 
the space of c\`adl\`ag paths in $\Omega$. We further write
\be{IPSlaw} 
P^\mu(\cdot) = \int_{\Omega} P^\eta(\cdot)\,\mu(\di\eta) 
\ee
to denote the law of $\xi$ when $\xi(0)$ is drawn from $\mu$. We further assume that
\be{tailtriv}
P^\mu \mbox{ is tail trivial}.
\ee

Conditional on $\xi$, let
\be{rwdef} 
X = (X_t)_{t\geq 0}
\ee
be the random walk with local transition rates
\be{rwtrans}
\begin{aligned}
&x \to x+1 \quad \mbox{ at rate } \quad \alpha\,\xi_t(x) + \beta\,[1-\xi_t(x)],\\
&x \to x-1 \quad \mbox{ at rate } \quad \beta\,\xi_t(x) + \alpha\,[1-\xi_t(x)],
\end{aligned}
\ee
where w.l.o.g. 
\be{albe}
0 < \beta < \alpha < \infty.
\ee 
In words, on occupied sites the random walk jumps to the right at rate $\alpha$ and 
to the left at rate $\beta$, while at vacant sites it does the opposite. Note that,
by (\ref{albe}), on occupied sites the drift is positive, while on vacant sites it is 
negative. Also note that the sum of the jump rates is $\alpha+\beta$ and is independent 
of $\xi$. For $x\in\Z$, we write $P^\xi_0$ to denote the law of $X$ starting from $X(0)=0$ 
conditional on $\xi$, and
\be{Pxnudef} 
\P_{\mu,0}(\cdot) = \int_{D_\Omega[0,\infty)}
P^\xi_0(\cdot)\,P^\mu(\di\xi)
\ee
to denote the law of $X$ averaged over $\xi$. We refer to $P^\xi_0$ as the \emph{quenched} 
law and to $\P_{\mu,0}$ as the \emph{annealed} law.


\subsection{Large deviation principles}
\label{S1.2} 

In \cite{AvdHoRe09} we proved that if $\xi$ is cone-mixing, then $X$ satisfies a law of large 
numbers (LLN), i.e., there exists a $v\in\R$ such that 
\be{LLN} 
\lim_{t\to\infty} t^{-1}X_t = v \qquad \P_{\mu,0}-a.s.
\ee
All attractive spin-flip systems for which the coupling time at the origin, starting from 
the configurations $\eta \equiv 1$ and $\eta \equiv 0$, has finite mean are cone-mixing. 
Theorems~\ref{aLDP}--\ref{qLDP} below state that $X$ satisfies both an annealed and a 
quenched large deviation principle (LDP); the interval $K$ in (\ref{Iannrate}) and 
(\ref{Iquerate}) can be either open, closed or half open and half closed.

Define
\be{Mepsdef}
\begin{aligned}
M &= \sum_{x \neq 0}\,\sup_{\eta\in\Omega}
|c(0,\eta)-c(0,\eta^x)|,\\
\epsilon &= \inf_{\eta\in\Omega} |c(0,\eta) + c(0,\eta^0)|.
\end{aligned}
\ee
The interpretation of (\ref{Mepsdef}) is that $M$ is a measure for the \emph{maximal} 
dependence of the transition rates on the states of single sites, while $\epsilon$ is 
a measure for the \emph{minimal} rate at which the states of single sites change. See 
\cite{Li85}, Section I.4, for examples. In \cite{AvdHoRe09} we showed that if $M<\epsilon$ 
then $\xi$ is cone-mixing.  

\bt{aLDP} {\bf [Annealed LDP]} Assume {\rm (\ref{shift}--\ref{erg})}.\\
(a) There exists a convex rate function $I^\mathrm{ann}\colon\,\R \to [0,\infty)$, 
satisfying
\be{Iann}
I^\mathrm{ann}(\theta)\,\,\left\{\begin{array}{ll}
= 0,  &\mbox{ if } \theta \in [v_-^\mathrm{ann},v_+^\mathrm{ann}],\\
> 0,  &\mbox{ if } \theta \in \R \backslash [v_-^\mathrm{ann},v_+^\mathrm{ann}],
\end{array}
\right. 
\ee 
for some $-(\alpha-\beta) \leq v_-^\mathrm{ann} \leq v \leq v_+^\mathrm{ann} \leq 
\alpha-\beta$, such that
\be{Iannrate}
\lim_{t\to\infty} \frac{1}{t} \log \P_{\mu,0}\big(t^{-1}X_t \in K\big)
= - \inf_{\theta \in K} I^\mathrm{ann}(\theta)
\ee
for all intervals $K$ such that either $K \varsubsetneq [v_-^\mathrm{ann},v_+^\mathrm{ann}]$ or 
$\mathrm{int}(K) \ni v$.\\  
(b) $\lim_{|\theta|\to\infty} I^\mathrm{ann}(\theta)/|\theta|=\infty$.\\ 
(c) If $M<\epsilon$ and $\alpha-\beta<\tfrac12(\epsilon-M)$, then
\be{Ann!0}
v_-^\mathrm{ann} = v = v_+^\mathrm{ann}.
\ee
\et

\bt{qLDP} {\bf [Quenched LDP]} Assume {\rm (\ref{shift}--\ref{erg}) and (\ref{tailtriv})}.\\
(a) There exists a convex rate function $I^\mathrm{que}\colon\,\R \to [0,\infty)$, 
satisfying
\be{Ique}
I^\mathrm{que}(\theta)\,\,\left\{\begin{array}{ll}
= 0,  &\mbox{ if } \theta \in [v_-^\mathrm{que},v_+^\mathrm{que}],\\
> 0,  &\mbox{ if } \theta \in \R \backslash [v_-^\mathrm{que},v_+^\mathrm{que}],
\end{array}
\right.
\ee 
for some $-(\alpha-\beta) \leq v_-^\mathrm{que} \leq v \leq v_+^\mathrm{que} \leq 
\alpha-\beta$, such that
\be{Iquerate}
\lim_{t\to\infty} \frac{1}{t} \log P^\xi\big(t^{-1}X_t \in K\big)
= - \inf_{\theta \in K} I^\mathrm{que}(\theta) \qquad \xi-a.s.
\ee
for all intervals $K$.\\ 
(b) $\lim_{|\theta|\to\infty} I^\mathrm{que}(\theta)/|\theta|=\infty$ and
\be{Iquesym} 
I^\mathrm{que}(-\theta) = I^\mathrm{que}(\theta)
+ \theta(2\rho-1)\log(\alpha/\beta), \qquad \theta \geq 0.
\ee
(c) If $M<\epsilon$ and $\alpha-\beta<\tfrac12(\epsilon-M)$, then
\be{Que!0}
v_-^\mathrm{que} = v = v_+^\mathrm{que}.
\ee
\et

\noindent
Theorems~\ref{aLDP} and \ref{qLDP} are proved in Sections~\ref{S2} and \ref{S3}, respectively. 
We are not able to show that (\ref{Iannrate}) holds for all closed intervals $K$, although we 
expect this to be true in general. 

Because
\be{Iquanrel}
I^\mathrm{que} \geq I^\mathrm{ann},
\ee
Theorems~\ref{qLDP}(b--c) follow from Theorems~\ref{aLDP}(b--c), with the exception of 
the symmetry relation (\ref{Iquesym}). There is no symmetry relation analogous to 
(\ref{Iquesym}) for $I^\mathrm{ann}$. It follows from (\ref{Iquanrel}) that
\be{vrels}
 v_-^\mathrm{ann} \leq v_-^\mathrm{que} \leq v \leq v_+^\mathrm{que} \leq 
v_+^\mathrm{ann}.
\ee


\subsection{Random walk in dynamic random environment: simple symmetric exclusion}
\label{S1.3}

It is natural to ask whether in a dynamic random environment the rate functions
always have a unique zero. The answer is no. In this section we show that when
$\xi$ is the \emph{simple symmetric exclusion process} in equilibrium with an 
arbitrary density of occupied sites $\rho \in (0,1)$, then for any $0<\beta<\alpha
<\infty$ the probability that $X_t$ is near the origin decays slower than exponential 
in $t$. Thus, slow-down is possible not only in a static random environment (see
Section \ref{S1.4}), but also in a dynamic random environment, provided it is not 
fast mixing. Indeed, the simple symmetric exclusion process is not even cone-mixing.
   
The one-dimensional simple symmetric exclusion process
\be{SSE}
\xi = \{\xi_t(x)\colon\,x\in\Z,\,t\geq 0\}
\ee
is the Markov process on state space $\Omega=\{0,1\}^\Z$ with generator $L$ given by
\be{ssegen}
(Lf)(\eta) = \sum_{ {x,y\in\Z} \atop {x \sim y} } [f(\eta^{xy})-f(\eta)],
\qquad \eta\in\Omega,
\ee
where $f$ is any cylinder function on $\R$, the sum runs over unordered neighboring pairs 
of sites in $\Z$, and $\eta^{xy}$ is the configuration obtained from $\eta$ by interchanging
the states at sites $x$ and $y$. We will asume that $\xi$ starts from the Bernoulli product 
measure with density $\rho \in (0,1)$, i.e., at time $t=0$ each site is occupied with 
probability $\rho$ and vacant with probability $1-\rho$. This measure, which we denote by 
$\nu_\rho$, is an equilibrium for the dynamics (see \cite{Li85}, Theorem VIII.1.44). 

Conditional on $\xi$, the random walk 
\be{RWonSSE}
X = (X_t)_{t\geq 0}
\ee
has the same local transition rates as in (\ref{rwtrans}--\ref{albe}). We also retain the
definition of the quenched law $P_0^\xi$ and the annealed law $\P_{\nu_\rho,0}$, as in
(\ref{Pxnudef}) with $\mu=\nu_\rho$.

Since the simple symmetric exclusion process is \emph{not} cone-mixing (the space-time 
mixing property assumed in \cite{AvdHoRe09}), we do not have the LLN. Since it is \emph{not} 
an attractive spin-flip system either, we also do not have the LDP. We plan to address 
these issues in future work. Our main result here is the following.

\bt{SSEslow}
For all $\rho \in (0,1)$,
\be{zerospeedrate}
\lim_{t\to\infty} \frac{1}{t} \log \P_{\nu_\rho,0}\big(|X_t| \leq 2\sqrt{t\log t}\big) = 0.
\ee
\et

\noindent
Theorem~\ref{SSEslow} is proved in Section~\ref{S4}.


\subsection{Discussion}
\label{S1.4} 

{\bf Literature.}
Random walk in \emph{static} random environment has been an intensive research area 
since the early 1970's. One-dimensional models are well understood. In particular,
recurrence vs.\ transience criteria, laws of large numbers and central theorems have 
been derived, as well as quenched and annealed large deviation principles. In higher 
dimensions a lot is known as well, but some important questions still remain open. For 
an overview of these results, we refer the reader to \cite{Ze01,Ze06} and \cite{Sz02}. 
See the homepage of Firas Rassoul-Agha [{\sc www.math.utah.edu/$\sim$firas/Research}] 
for an up-to-date list of references. 

For random walk in \emph{dynamic} random environment the state of the art is rather more 
modest, even in one dimension. Early work was done in \cite{Ma86}, which considers a 
one-dimensional environment consisting of spins flipping independently between $-1$ and 
$+1$, and a walk that at integer times jumps left or right according to the spin it sees 
at that time. A necessary and sufficient criterion for recurrence is derived, as well as 
a law of large numbers.

Three classes of models have been studied in the literature so far:
\begin{itemize}
\item[(1)] 
\emph{Space-time random environment}: globally updated at each unit of time 
\cite{BoMiPe97,BoMiPe00,BoMiPe04,Be04,RaAgSe05,Yi09};
\item[(2)] 
\emph{Markovian random environment}: independent in space and locally updated according to 
a single-site Markov chain 
\cite{BoMiPe00,Ic02,BaZe06}; 
\item[(3)]
\emph{Weak random environment}: small perturbation of homogeneous random walk (possibly with 
a feedback of the walk on the environment)
\cite{BoMiPe94a,BoMiPe94b,BoMiPe94c,BeBoMiPe96,BeBoMiPe98,IgRo98,BoMiPe07}.
\end{itemize}
The focus of these references is: transience vs.\ recurrence \cite{Ma86,Ic02}, central limit 
theorem \cite{BoMiPe94a,BeBoMiPe96,BoMiPe97,BeBoMiPe98,BoMiPe00,BoMiPe04,Be04,RaAgSe05,BoMiPe07},
law of large numbers and central limit theorem \cite{BaZe06}, decay of correlations in 
space and time \cite{BoMiPe94b}, convergence of the law of the environment as seen from 
the walk \cite{BoMiPe94c}, large deviations \cite{IgRo98,Yi09}. In classes (1) and (2) 
the random environment is uncorrelated in time, respectively, in space. In \cite{AvdHoRe09} we 
moved away from this restriction by proving a law of large numbers for a class of dynamic 
random environments correlated in space \emph{and} time, satisfying a space-time mixing condition 
called \emph{cone-mixing}. We showed that a large class of uniquely ergodic attractive spin-flip 
systems falls into this class.

Consider a \emph{static} random environment $\eta$ with law $\nu_\rho$, the Bernoulli product 
measure with density $\rho \in (0,1)$, and a random walk $X=(X_t)_{t \geq 0}$ with transition
rates (compare with (\ref{rwtrans})) 
\be{rwtransdis}
\begin{aligned}
&x \to x+1 \quad \mbox{ at rate } \quad \alpha\eta(x)+\beta[1-\eta(x)],\\
&x \to x-1 \quad \mbox{ at rate } \quad \beta\eta(x)+\alpha[1-\eta(x)],
\end{aligned}
\ee
where $0<\beta<\alpha<\infty$. In \cite{So75} it is shown that $X$ is recurrent when $\rho
=\tfrac12$ and transient to the right when $\rho>\tfrac12$. In the transient case both 
ballistic and non-ballistic behavior occur, i.e., $\lim_{t\to\infty} X_t/t=v$ for 
$\P_{\nu_\rho}$-a.e.\ $\xi$, and
\be{statspeed}
v\,\,\left\{\begin{array}{ll}
= 0 &\mbox{ if } \rho \in [\frac12,\rho_c],\\
> 0 &\mbox{ if } \rho \in (\rho_c,1],
\end{array}
\right.
\ee
where
\be{rhocdef} 
\rho_c = \frac{\alpha}{\alpha+\beta} \in (\tfrac12,1), 
\ee
and, for $\rho \in (\rho_c,1]$,
\be{vcrhoc}
v = v(\rho,\alpha,\beta) =
(\alpha+\beta)\,\frac{\alpha\beta+\rho(\alpha^2-\beta^2)-\alpha^2}
{\alpha\beta-\rho(\alpha^2-\beta^2)+\alpha^2} =
(\alpha-\beta)\,\frac{\rho-\rho_c}{\rho(1-\rho_c)+\rho_c(1-\rho)}.
\ee

\medskip\noindent
{\bf Attractive spin flips.}
The analogues of (\ref{Iannrate}) and (\ref{Iquerate}) in the \emph{static} random environment 
(with no restriction on the interval $K$ in the annealed case) were proved in \cite{GedHo94} 
(quenched) and \cite{CoGaZe00} (quenched and annealed). Both $I^\mathrm{ann}$ and $I^\mathrm{que}$ 
are zero on the interval $[0,v]$ and are strictly positive outside (``slow-down phenomenon''). 
For $I^\mathrm{que}$ the same symmetry property as in (\ref{Iquesym}) holds. Moreover, an 
explicit formula for $I^\mathrm{que}$ is known in terms of random continued fractions. 

We do not have explicit expressions for $I^\mathrm{ann}$ and $I^\mathrm{que}$ in the 
\emph{dynamic} random environment. Even the characterization of their zero sets remains 
open, although under the stronger assumptions that $M<\epsilon$ and $\alpha-\beta<\tfrac12
(\epsilon-M)$ we know that both have a unique zero at $v$.

Theorems~\ref{aLDP}--\ref{qLDP} can be generalized beyond spin-flip systems, i.e., systems 
where more than one site can flip state at a time. We will see in Sections~\ref{S2}--\ref{S3} 
that what really matters is that the system has positive correlations in space and time. As 
shown in \cite{Ha77}, this holds for monotone systems (see \cite{Li85}, Definition II.2.3) if 
and only if all transitions are such that they make the configuration either larger or smaller 
in the partial order induced by inclusion.

\medskip\noindent
{\bf Simple symmetric exclusion.}
What Theorem~\ref{SSEslow} says is that, for all choices of the parameters, the annealed 
rate function (if it exists) is zero at $0$, and so there is a slow-down phenomenon similar 
to what happens in the static random environment. We will see in Section~\ref{S4} that this 
slow-down comes from the fact that the simple symmetric exclusion process suffers ``traffic 
jams'', i.e., long strings of occupied and vacant sites have an appreciable probability to 
survive for a long time.

\begin{figure}[hbtp]
\vspace{0.5cm}
\begin{center}
\includegraphics[width=6cm]{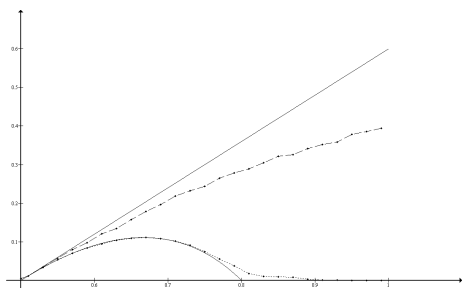}
\end{center}
\caption{\small Speeds as a function of $p$ for $\rho=0.8$.} 
\label{fig-rho}
\end{figure}

\begin{figure}[hbtp]
\vspace{0.1cm}
\begin{center}
\includegraphics[width=6cm]{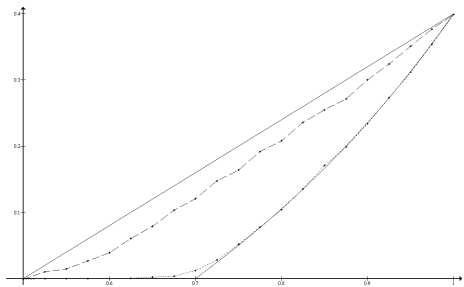}
\end{center}
\caption{\small Speeds as a function of $\rho$ for $p=0.7$.}
\label{fig-p} 
\vspace{0.3cm}
\end{figure}

To test the validity of the LLN for the simple symmetric exclusion process, we performed a 
simulation the outcome of which is drawn in Figs.~\ref{fig-rho}--\ref{fig-p}. For each point 
in these figures, we drew $10^3$ initial configurations according to the Bernoulli 
product measure with density $\rho$, and from each of these configurations ran a 
discrete-time exclusion process with parallel updating for $10^4$ steps. Given the 
latter, we ran a discrete-time random walk for $10^4$ steps, both in the static 
environment (ignoring the updating) and in the dynamic environment (respecting the 
updating), and afterwards averaged the displacement of the walk over the $10^3$ 
initial configurations. The probability to jump to the right was taken to be 
$p$ on an occupied site and $q=1-p$ on a vacant site, where $p$ replaces $\alpha/
(\alpha+\beta)$ in the continuous-time model. In Figs.~\ref{fig-rho}--\ref{fig-p}, 
the speeds resulting from these simulations are plotted as a function of $p$ for 
$\rho=0.8$, respectively, as a function of $\rho$ for $p=0.7$. In each figure we 
plot four curves: (1) the theoretical speed in the static case (as described by 
(\ref{vcrhoc})); (2) the simulated speed in the static case; (3) the simulated 
speed in the dynamic case; (4) the speed for the average environment, i.e., 
$(2\rho-1)(2p-1)$. The order in which these curves appear in the figures is 
from bottom to top.


Fig.~\ref{fig-rho} shows that, in the static case with $\rho$ fixed, as $p$ increases 
the speed first goes up (because there are more occupied than vacant sites), and then 
goes down (because the vancant sites become more efficient to act as a barrier). In the 
dynamic case, however, the speed is an increasing function of $p$: the vacant site are 
not frozen but move around and make way for the walk. It is clear from Fig.~\ref{fig-p} 
that the only value of $\rho$ for which there is a zero speed in the dynamic case is 
$\rho=\tfrac12$, for which the random walk is recurrent. Thus, the simulation suggests
that there is no (!) non-ballistic behavior in the transient case. In view of 
Theorem~\ref{SSEslow}, this in turn suggests that the annealed rate function (if it 
exists) has zero set $[0,v]$. 

In both pictures the two curves at the bottom should coincide. Indeed, they almost 
coincide, except for values of the parameters that are close to the transition 
between ballistic and non-ballistic behavior, for which fluctuations are to be 
expected. Note that the simulated speed in the dynamic environment lies inbetween 
the speed for the static environment and the speed for the average environment. 
We may think of the latter as corresponding to a simple symmetric exclusion process
running at rate $0$, respectively, $\infty$ rather than at rate $1$ as in (\ref{ssegen}).


\section{Proof of Theorem~\ref{aLDP}}
\label{S2}

In Section \ref{S2.1} we prove three lemmas for the probability that the empirical speed 
is above a given threshold. These lemmas will be used in Section \ref{S2.2} to prove 
Theorems~\ref{aLDP}(a--b). In Section \ref{S2.3} we prove Theorems~\ref{aLDP}(c).


\subsection{Three lemmas}
\label{S2.1}

\bl{Psubad} For all $\theta\in\R$,
\be{Jsubad}
J^+(\theta) = -\lim_{t\to\infty} \frac{1}{t} \log \P_{\mu,0}(X_t \geq \theta t)
\quad \mbox{exist and is finite}.
\ee 
\el

\begin{proof}
For $z\in\Z$ and $u \geq 0$, let $\sigma_{z,u}$ denote the operator acting on $\xi$ as 
\be{sigtdef} 
(\sigma_{z,u}\xi)(x,t) = \xi(z+x,u+t), \qquad x\in\Z,\,t\geq 0. 
\ee
Fix $\theta \neq 0$, and let $G_\theta=\{t \geq 0\colon\,\theta t\in\Z\}$ be the non-negative
grid of width $1/|\theta|$. For any $s,t \in G_\theta$, we have
\be{Psubad3}
\begin{aligned}
&\P_{\mu,0}\big(X_{s+t} \geq \theta(s+t)\big)
= E^\mu\Big[P_0^\xi\big(X_{s+t} \geq \theta(s+t)\big)\Big]\\[0.2cm]
&= \sum_{y\in\Z}E^\mu\Big[P_0^\xi(X_s=y)\,
P_y^{\sigma_{0,s}\xi}\big(X_t\geq \theta (s+t)\big)\Big] \geq
\sum_{y\geq\theta s}E^\mu\Big[P_0^\xi(X_s=y)\,
P_{\theta s}^{\sigma_{0,s}\xi}\big(X_t\geq\theta(s+t)\big)\Big]\\[0.2cm]
&= E^\mu\Big[P_0^\xi(X_s\geq\theta s)\, P_0^{\sigma_{\theta s,s}\xi}
\big(X_t\geq\theta t\big)\Big] \geq E^\mu\Big[P_0^\xi(X_s\geq\theta s)\Big]\,
E^\mu\Big[P_0^{\sigma_{\theta s,s}\xi}\big(X_t\geq\theta t\big)\Big]\\[0.2cm]
&=\P_{\mu,0}(X_s\geq\theta s)\,\P_{\mu,0}(X_t\geq\theta t).
\end{aligned}
\ee
The first inequality holds because two copies of the random walk running on the same realization 
of the random environment can be coupled so that they remain ordered. The second inequality uses 
that
\be{Psubad5} 
\xi\mapsto P_0^\xi(X_s\geq\theta s) \quad \mbox{and} \quad 
\xi\mapsto P_0^{\sigma_{\theta s,s}\xi}\big(X_t\geq\theta t\big) 
\ee 
are non-decreasing and that the law $P^\mu$ of an attractive spin-flip system has the FKG-property 
in space-time (see \cite{Li85}, Corollary II.2.12). Let
\be{gdefsunad} 
g(t) = -\log \P_{\mu,0}(X_t\geq\theta t).
\ee
Then it follows from (\ref{Psubad3}) that $(g(t))_{t\geq 0}$ is subadditive along $G_\theta$, i.e.,
$g(s+t) \leq g(s)+g(t)$ for all $s,t\in G_\theta$. Since $\P_{\mu,0}(X_t \geq \theta t)>0$ for 
all $t\geq 0$, it thefore follows that
\be{Jsubadgrid}
J^+(\theta) = -\lim_{ {t\to\infty} \atop {t \in G_\theta} } 
\frac{1}{t} \log \P_{\mu,0}(X_t \geq \theta t) \quad \mbox{exist and is finite}.
\ee
Because $X$ takes values in $\Z$, the restriction $t \in G_\theta$ can be removed. This proves 
the claim for $\theta \neq 0$. The claim easily extends to $\theta=0$, because the transition 
rates of the random walk are bounded away from $0$ and $\infty$ uniformly in $\xi$ (recall
(\ref{rwtrans})). 
\end{proof}

\bl{Jprops} 
$\theta \mapsto J^+(\theta)$ is non-decreasing and convex on $\R$. 
\el

\begin{proof}
We follow an argument similar to that in the proof of Proposition~\ref{Psubad}. Fix $\theta,
\gamma\in\R$ and $p\in [0,1]$ such that $p\gamma,(1-p)\theta\in\Z$. Estimate
\be{J2}
\begin{aligned}
&\P_{\mu,0}\big(X_t \geq [p\gamma+(1-p)\theta] t\big)
= E^\mu\Big[P_{0}^\xi\big(X_t \geq [p\gamma+(1-p)\theta] t\big)\Big]\\[0.2cm]
&\qquad = \sum_{y\in\Z} E^\mu\Big[P_0^\xi(X_{pt}=y)\,P_y^{\sigma_{0,pt}\xi}
\big(X_{t(1-p)} \geq [p\gamma+(1-p)\theta] t\big)\Big]\\[0.2cm]
&\qquad \geq \sum_{y\geq p\gamma t}
E^\mu\Big[P_0^\xi(X_{pt}=y)\, P_{p\gamma t}^{\sigma_{0,pt}\xi}
\big(X_{t(1-p)} \geq [p\gamma+(1-p)\theta] t\big)\big]\\[0.2cm]
&\qquad =E^\mu\Big[P_0^\xi(X_{pt}\geq p\gamma t)\,
P_0^{\sigma_{p\gamma t,pt}\xi}
\big(X_{t(1-p)} \geq (1-p)\theta t\big)\Big]\\[0.2cm]
&\qquad\geq E^\mu\Big[P_0^\xi(X_{pt}\geq p\gamma t)\Big]\,
E^\mu\Big[P_{p\gamma\,t}^{\sigma_{p\gamma t,pt}\xi}
\big(X_{t(1-p)}\geq (1-p)\theta t\big)\Big]\\[0.2cm]
&\qquad = \P_{\mu,0}(X_{pt} \geq p\gamma t)\,
\P_{\mu,0}\big(X_{t(1-p)} \geq (1-p)\theta t\big).
\end{aligned}
\ee
It follows from (\ref{J2}) and the remark below (\ref{Jsubadgrid}) that 
\be{J4} 
-J^+\big(p\gamma+(1-p)\theta\big) \geq -pJ^+(\gamma)-(1-p)J^+(\theta), 
\ee 
which settles the convexity.
\end{proof}

\bl{Jnonzero} 
$J^+(\theta)>0$ for $\theta>\alpha-\beta$ and $\lim_{\theta\to\infty} J^+(\theta)/\theta
=\infty$.
\el

\begin{proof}
Let $(Y_t)_{t\geq 0}$ be the nearest-neighbor random walk on $\Z$ that jumps to the right at 
rate $\alpha$ and to the left at rate $\beta$. Write $\P^{\RW}_0$ to denote its law starting 
from $Y(0)=0$. 
Clearly,
\be{I1} 
\P_{\mu,0}(X_t \geq \theta t) \leq \P^{\RW}_0(Y_t \geq \theta t) \qquad \forall\,\theta\in\R.
\ee
Moreover,
\be{I2} 
J^{\RW}(\theta) = -\lim_{t\to\infty} \frac{1}{t}
\log\P^{\RW}_0(Y_t \geq \theta t) 
\ee
exists, is finite and satisfies
\be{I3} 
J^{\RW}\left(\alpha-\beta\right) = 0, 
\quad
J^{\RW}(\theta)>0 \mbox{ for }\theta>\alpha-\beta,
\quad 
\quad \lim_{\theta\to\infty} J^{\RW}(\theta)/\theta=\infty. 
\ee
Combining (\ref{I1}--\ref{I3}), we get the claim. 
\end{proof}

Lemmas~\ref{Psubad}--\ref{Jnonzero} imply that an \emph{upward} annealed LPD holds with 
a rate function $J^+$ whose qualitative shape is given in Fig.~\ref{fig-J^+}. 


\begin{figure}[hbtp]
\vspace{1cm}
\begin{center}
\setlength{\unitlength}{0.3cm}
\begin{picture}(12,10)(-3,-0.8)
\put(-12,0){\line(24,0){24}}
\put(0,0){\line(0,10){10}}
{\thicklines
\qbezier(-8,0)(0,0)(5,0)
\qbezier(5,0)(7,0.2)(10,7) }
\put(-.3,-1.5){$v$} 
\put(4,-1.5){$v_+^\mathrm{ann}$}
\put(12.8,0){$\theta$} 
\put(-1,10.8){$J^+(\theta)$}
\put(5,0){\circle*{.5}}
\end{picture}
\end{center}
\caption{\small Shape of $\theta \mapsto J^+(\theta)$.}
\label{fig-J^+}
\end{figure}



\begin{figure}[hbtp]
\vspace{1cm}
\begin{center}
\setlength{\unitlength}{0.3cm}
\begin{picture}(12,10)(-3,-0.8)
\put(-12,0){\line(24,0){24}} 
\put(0,0){\line(0,10){10}}
{\thicklines 
\qbezier(8,0)(0,0)(-5,0) 
\qbezier(-5,0)(-7,0.2)(-10,7)
}
\put(-.3,-1.5){$v$} 
\put(-5.7,-1.5){$v_-^\mathrm{ann}$}
\put(12.8,0){$\theta$} 
\put(-1,10.8){$J^-(\theta)$}
\put(-5,0){\circle*{.5}}
\end{picture}
\end{center}
\caption{\small Shape of $\theta \mapsto J^-(\theta)$.}
\label{fig-J^-}
\end{figure}


\subsection{Annealed LDP}
\label{S2.2}

Clearly, $J^+$ depends on $P^\mu$, $\alpha$ and $\beta$. Write 
\be{J+}
J^+=J_{P^\mu,\alpha,\beta}
\ee
to exhibit this dependence. So far we have not used the restriction $\alpha>\beta$ in 
(\ref{albe}). By noting that $-X_t$ is equal in distribution to $X_t$ when $\alpha$ 
and $\beta$ are swapped and $P^\mu$ is replaced by $\bar{P}^\mu$, the image of $P^\mu$
under reflection in the origin (recall (\ref{rwtrans})), we see that the \emph{upward} 
annealed LDP proved in Section~\ref{S2.1} also yields a \emph{downward} annealed LDP 
\be{J-lim}
J^-(\theta) = - \lim_{t\to\infty} \frac{1}{t} \log \P_{\mu,0}(X_t \leq \theta t),
\qquad \theta \in \R,
\ee
with 
\be{J-}
J^-=J_{\bar{P}^\mu,\beta,\alpha},
\ee 
whose qualitative shape is given in Fig.~\ref{fig-J^-}. Note that
\be{vv-v+rel}
v_-^\mathrm{ann} \leq v \leq  v_+^\mathrm{ann},
\ee
because $v$, the speed in the LLN proved in \cite{AvdHoRe09}, must lie in the zero set 
of both $J^+$ and $J^-$. 

Our task is to turn the upward and downward annealed LDP's into the annealed LDP of 
Theorem~\ref{aLDP}. 

\bp{AnnLDP} 
Let
\be{Iann*}
I^\mathrm{ann}(\theta) = \left\{\begin{array}{ll}
J_{P^\mu,\alpha,\beta}(\theta)   &\mbox{ if } \theta \geq v,\\
J_{\bar{P}^\mu,\beta,\alpha}(-\theta)   &\mbox{ if } \theta \leq v.
\end{array}
\right. 
\ee
Then
\be{Iannratealt}
\lim_{t\to\infty} \frac{1}{t} \log \P_{\mu,0}\big(t^{-1}X_t \in K)
= - \inf_{\theta \in K} I^\mathrm{ann}(\theta)
\ee
for all closed intervals such that either $K\varsubsetneq [v_-^\mathrm{ann},v_+^\mathrm{ann}]$ or
$\mathrm{int}(K) \ni v$.
\ep

\begin{proof}
We distinguish three cases.

\medskip\noindent
(1) $K \subset [v,\infty)$, $K \varsubsetneq [v, v_+^\mathrm{ann}]$: Let $\mathrm{cl}(K)=[a,b]$.
Then, because $J^+$ is continuous,
\be{I1*}
\frac{1}{t} \log \P_{\mu,0}\big(t^{-1}X_t \in K\big)
= \frac{1}{t} \log \Big[e^{-tJ^+(a)+o(t)}- e^{-tJ^+(b)+o(t)}\Big].
\ee
By Lemma \ref{Jprops}, $J^+$ is strictly increasing on $[v_+^\mathrm{ann},\infty)$, and so 
$J^+(b)>J^+(a)$. Letting $t\to\infty$ in (\ref{I1*}), we therefore see that 
\be{case1} 
\lim_{t\to\infty} \frac{1}{t} \log \P_{\mu,0}\big(t^{-1}X_t \in K\big) 
= -J^+(a) = -\inf_{\theta \in K} I^\mathrm{ann}(\theta).
\ee

\medskip\noindent
(2) $K \subset (-\infty,v]$, $K \varsubsetneq [v_-^\mathrm{ann},v]$: Same as for (1) 
with $J^-$ replacing $J^+$.

\medskip\noindent
(3) $\mathrm{int}(K) \ni v$: In this case (\ref{Iannratealt}) is an immediate consequence of the 
LLN in (\ref{LLN}).
\end{proof}

\noindent
Proposition~\ref{AnnLDP} completes the proof of Theorems~\ref{aLDP}(a--b). Recall
(\ref{J+}) and (\ref{J-}). The restriction on $K$ comes from the fact that the difference
of two terms that are both $\exp[o(t)]$ may itself not be $\exp[o(t)]$.  


\subsection{Unique zero of $I^\mathrm{ann}$ when $M<\eps$}
\label{S2.3}

In \cite{AvdHoRe09} we showed that if $M<\epsilon$ and $\alpha-\beta<\tfrac12(\eps-M)$,
then a proof of the LLN can be given that is based on a perturbation argument for the
generator of the \emph{environment process}
\be{envproc}
\zeta = (\zeta_t)_{t \geq 0}, \qquad \zeta_t = \tau_{X_t}\xi_t, 
\ee
i.e., the random environment as seen relative to the random walk. In particular, it is 
shown that $\zeta$ is uniquely ergodic with equilibrium $\mu_e$. This leads to a series 
expansion for $v$ in powers of $\alpha-\beta$, with coefficients that are functions of 
$P^\mu$ and $\alpha+\beta$ and that are computable via a recursive scheme. The speed in 
the LLN is given by 
\be{vpert}
v = (2\widetilde{\rho}-1)(\alpha-\beta)
\ee
with $\widetilde{\rho} = \langle\eta(0)\rangle_{\mu_e}$, where $\langle\cdot
\rangle_{\mu_e}$ denotes expectation over $\mu_e$ ($\widetilde{\rho}$ is the 
fraction of time $X$ spends on occupied sites).

\bp{Unique0} 
Let $\xi$ be an attractive spin-flip system with $M<\epsilon$. If $\alpha-\beta<\tfrac12
(\epsilon-M)$, then the rate function $I^\mathrm{ann}$ in {\rm (\ref{I1*})} has a unique 
zero at $v$. 
\ep

\begin{proof}
It suffices to show that 
\be{LDexp}
\limsup_{t\to\infty} \frac{1}{t} \log \P_{\mu,0}\big(|t^{-1}X_t-v| \geq 2\delta\big) < 0
\qquad \forall\,\delta>0. 
\ee
To that end, put $\gamma=\delta/2(\alpha-\beta)>0$. Then, by (\ref{vpert}), $v\pm\delta
=[2(\widetilde{\rho}\pm\gamma)-1](\alpha-\beta)$. Let
\be{OccTime} 
A_t = \int_0^t \xi_s(X_s)\,\di s
\ee
be the time $X$ spends on occupied sites up to time $t$, and define
\be{E}
E_t= \big\{|t^{-1}A_t-\widetilde{\rho}| \geq \gamma\big\}.
\ee
Estimate
\be{DeltaConv}
\P_{\mu,0}(|t^{-1}X_t-v| \geq 2\delta\big)
\leq \P_{\mu,0}(E_t) + \P_{\mu,0}\big(|t^{-1}X_t-v| \geq 2\delta \mid E_t^c\big).
\ee
Conditional on $E_t^c$, $X$ behaves like a homogeneous random walk with speed in 
$[v-\delta,v+\delta]$. Therefore the second term in the r.h.s.\ of (\ref{DeltaConv}) 
vanishes exponentially fast in $t$. In \cite{AvdHoRe09}, Lemma 3.4, Eq.\ (3.26) and 
Eq.\ (3.36), we proved that
\be{Bounds}
\trn S(t)f \trn \leq \erom^{-c_1 t}\,\trn f \trn
\quad \text{ and } \quad 
\left\| S(t)f - \left\langle f \right\rangle_{\mu_e}\right\|_\infty 
\leq c_2\,\erom^{-(\epsilon-M)t}\,\trn f \trn
\ee
for some $c_1,c_2 \in (0,\infty)$, where $S=(S(t))_{t\geq 0}$ denotes the semigroup
associated with the environment process $\zeta$, and $\trn f \trn$ denotes the triple 
norm of $f$. As shown in \cite{ReVo09}, (\ref{Bounds}) implies a Gaussian concentration 
bound for additive functionals, namely, 
\be{concbd}
\P_{\mu,0}\left(\left| t^{-1} \int_0^t f(\zeta_s) 
- \langle f \rangle_\mu\right| \geq \gamma\right)
\leq c_3\,\erom^{-\gamma^2 t/c_4\,\trn f\trn^2}
\ee
for some $c_3,c_4 \in (0,\infty)$, uniformly in $t>0$, $f$ with $\trn f\trn < \infty$ 
and $\gamma>0$. By picking $f(\eta) = \eta(0)$, $\eta\in\Omega$, we get
\be{BoundE} 
\P_{\mu,0}\big(E_t\big) \leq c_5\,\erom^{-c_6 t}
\ee
for some $c_5,c_6 \in (0,\infty)$. Therefore also the first term in the r.h.s.\ of 
(\ref{DeltaConv}) vanishes exponentially fast in $t$.
\end{proof}

\noindent
Proposition~\ref{Unique0} completes the proof of Theorems~\ref{aLDP}(c).


\section{Proof of Theorem~\ref{qLDP}}
\label{S3}

In Section~\ref{S3.1} we prove three lemmas for the probability that the empirical speed 
equals a given value. These lemmas will be used in Section \ref{S3.2} to prove 
Theorems~\ref{qLDP}(a--b). In Section \ref{S3.3} we prove Theorem~\ref{qLDP}(c).
Theorem~\ref{qLDP}(d) follows from Theorem~\ref{aLDP}(c) because $I^\mathrm{que}
\geq I^\mathrm{ann}$. 


\subsection{Three lemmas}
\label{S3.1}

In this section we state three lemmas that are the analogues of Lemmas 
\ref{Psubad}--\ref{Jnonzero}.

\bl{Psubad*} 
For all $\theta\in\R$, 
\be{Isubad*}
I^\mathrm{que}(\theta) = - \lim_{t\to\infty} \frac{1}{t} 
\log P_0^\xi(X_t = \lfloor\theta t\rfloor) 
\quad \mbox{exists, is finite and is constant $\xi$-a.s.}  
\ee 
\el

\begin{proof}
Fix $\theta \neq 0$, and recall that $G_\theta=\{t \geq 0\colon\,\theta t\in\Z\}$ is the 
non-negative grid of width $1/|\theta|$. For any $s,t \in G_\theta$, we have
\be{Psubad3*}
\begin{aligned}
P_0^\xi\big(X_{s+t} = \theta(s+t)\big) 
&\geq P_0^\xi\big(X_s = \theta s\big)\,
P_0^\xi\big(X_{s+t} = \theta(s+t) \mid X_s = \theta s\big)\\
&= P_0^\xi\big(X_s = \theta s\big)\,P_0^{T_s\xi}(X_t = \theta t),
\end{aligned}
\ee
where $T_s=\sigma_{\theta s,s}$. Let
\be{gdef}
g_t(\xi) = -\log P_0^{\xi}(X_t = \theta t).
\ee
Then it follows from (\ref{Psubad3*}) that $(g_t(\xi))_{t\geq 0}$ is a subadditive random 
process along $G_\theta$, i.e., $g_{s+t}(\xi) \leq g_s(\xi)+ g_t(T_s\xi)$ for all $s,t \in 
G_\theta$. From Kingman's subadditive ergodic theorem (see e.g.\ \cite{St89}) it therefore 
follows that
\be{Kingman}
\lim_{ {t\to\infty} \atop {t \in G_\theta} } \frac{1}{t}
\log P_0^{\xi}(X_t = \theta t) = - I^\mathrm{que}(\theta)
\ee
exists, is finite $\xi$-a.s, and is $T_s$-invariant for every $s \in G_\theta$. Moreover, 
since $\xi$ is ergodic under space-time shifts (recall (\ref{erg}) and (\ref{tailtriv})), 
this limit is constant $\xi$-a.s. Because the transition rates of the random walk are 
bounded away from $0$ and $\infty$ uniformly in $\xi$ (recall (\ref{rwtrans})), the 
restriction $t \in G_\theta$ may be removed after $X_t=\theta t$ is replaced by 
$X_t=\lfloor \theta t\rfloor$ in (\ref{Kingman}). This proves the claim for $\theta\neq 0$. 
By the boundedness of the transition rates, the claim easily extends to $\theta=0$.  
\end{proof}

\bl{Iprops*} 
$\theta \mapsto I^\mathrm{que}(\theta)$ is convex on $\R$. 
\el

\begin{proof}
The proof is similar to that of Proposition \ref{Psubad}. Fix $\theta,\zeta\in\R$ and 
$p \in [0,1]$. For any $t \geq 0$ such that $p\zeta t,(1-p)\theta t \in\Z$, we have
\be{J2*}
\begin{aligned}
P_0^\xi\big(X_t \geq [p\zeta+(1-p)\theta] t\big) &\geq
P_0^\xi\big(X_{pt} = p\zeta t\big)\,
P_0^\xi\big(X_t = [p\zeta+(1-p)\theta] t \mid X_{pt}=p\zeta t\big)\\
&= P_0^\xi\big(X_{pt} = p\zeta t\big)\, 
P_0^{\sigma_{p\zeta t,pt}\xi}\big(X_{(1-p)t} = (1-p)\theta t\big).
\end{aligned}
\ee
It follows from (\ref{J2*}) and the remark below (\ref{Jsubadgrid}) that
\be{J4*}
-I^\mathrm{que}\big(p\zeta+(1-p)\theta\big) 
\geq -pI^\mathrm{que}(\zeta) - (1-p)I^\mathrm{que}(\theta), 
\ee
which settles the convexity. 
\end{proof}

\bl{Inonzero} 
$I^\mathrm{que}(\theta)>0$ for $|\theta|>\alpha-\beta$ and $\lim_{\theta\to\infty}
I^\mathrm{que}(\theta)/|\theta|=\infty$. 
\el

\begin{proof}
Same as Lemma \ref{Jnonzero}.
\end{proof}


\subsection{Quenched LDP}
\label{S3.2}

We are now ready to prove the quenched LDP.

\bp{QLDP} For $P^\mu$-a.e.\ $\xi$, the family of probability measures $P^{\xi}_0
(X_t/t\in\,\cdot\,)$, $t>0$, satisfies the LDP with rate $t$ and with deterministic 
rate function $I^\mathrm{que}$.
\ep

\begin{proof}
Use Lemmas \ref{Psubad*}--\ref{Inonzero}.
\end{proof}

\noindent
Proposition~\ref{QLDP} completes the proof of Theorems~\ref{qLDP}, except for the 
symmetry relation in (\ref{Iquesym}), which will be proved in Section~\ref{S3.3}.
Recall (\ref{Iquanrel}) and the remark below it.


\subsection{A quenched symmetry relation}
\label{S3.3}

\bp{Symmetry} 
For all $\theta \in \R$, the rate function in Theorem {\rm \ref{QLDP}} satisfies 
the symmetry relation
\be{symmetry} 
I^\mathrm{que}(-\theta) = I^\mathrm{que}(\theta)+\theta(2\rho-1) \log (\alpha/\beta). 
\ee
\ep

\begin{proof}
We first consider a discrete-time random walk, i.e., a random walk that observes the
random environment and jumps at integer times. Afterwards we will extend the argument 
to the continuous-time random walk defined in (\ref{rwdef}--\ref{albe}). 

\medskip\noindent
{\bf 1.\ Path probabilities.}
Let 
\be{rwdefdis}
X=(X_n)_{n\in\N_0}
\ee 
be the random walk with transition probablities
\be{rwtransdiscr}
\begin{aligned}
&x \to x+1 \quad \mbox{ with probability } \quad p\,\xi_n(x) + q\,[1-\xi_n(x)],\\
&x \to x-1 \quad \mbox{ with probability } \quad q\,\xi_n(x) + p\,[1-\xi_n(x)],
\end{aligned}
\ee
where w.l.o.g.\ $p>q$. For an oriented edge $e=(i,i\pm1)$, $i\in \Z$, write $\e = (i \pm 1,i)$ 
to denote the reverse edge. Let $p_n(e)$ denote the probability for the walk to jump along the 
edge $e$ at time $n$. Note that in the \emph{static} random environment these probabilities are 
time-independent, i.e., $p_n(e)=p_0(e)$ for all $n\in\N$. 

We will be interested in $n$-step paths $\omega=(\omega_0,\ldots,\omega_n)\in\Z^n$ with 
$\omega_0=0$ and $\omega_n=\lfloor\theta n\rfloor$ for a given $\theta \neq 0$. Write 
$\Theta\omega$ to denote the time-reversed path, i.e., $\Theta\omega=(\omega_n,\ldots,\omega_0)$. 
Let $N_e(\omega)$ denote the number of times the edge $e$ is crossed by $\omega$, and write 
$t^j_e(\omega)$, $j=1,\ldots,N_e(\omega)$, to denote the successive times at which the edge 
$e$ is crossed. Let $E(\omega)$ denote the set of edges in the path $\omega$, and $E^+(\omega)$ 
the subset  of forward edges, i.e., edges of the form $(i,i+1)$. Then we have
\be{bombom} 
N_e (\Theta\omega) = N_{\e}(\omega) 
\ee
and
\be{bambam} 
t^j_e (\Theta\omega) = n+1-t^{N_{\e}(\omega)+1-j}_{\e}(\omega), 
\qquad j=1,\ldots,N_e(\Theta\omega) = N_{\e}(\omega). 
\ee
Given a realization of $\xi$, the probability that the walk follows the path $\omega$ equals
\be{propa} 
P^{\xi}(\omega) = \prod_{e\in E(\omega)} \prod_{j=1}^{N_e(\omega)} p_{t^j_e (\omega)}(e) 
= \prod_{e\in E^+(\omega)} \prod_{j=1}^{N_e(\omega)} p_{t^j_e (\omega)}(e)
\prod_{j=1}^{N_{\e}(\omega)} p_{t^j_{\e} (\omega)}(\e). 
\ee
The probability of the reversed path is, by (\ref{bombom}--\ref{bambam}),
\be{prore}
\begin{aligned}
P^{\xi}(\Theta\omega) 
&= \prod_{e\in E(\omega)} \prod_{j=1}^{N_e(\Theta\omega)} 
p_{t^j_e (\Theta\omega)}(e)
= \prod_{e\in E(\omega)} \prod_{j=1}^{N_{\e}(\omega)} 
p_{n+1- t_{\e}^{N_{\e}(\omega)+1-j}(\omega)}(e)\\
&= \prod_{e\in E(\omega)} \prod_{j=1}^{N_{\e}(\omega)} 
p_{n+1-t^j_{\e}(\omega)}(e) 
= \prod_{e\in E^+(\omega)} \prod_{j=1}^{N_e(\omega)} p_{n+1-t^j_e (\omega)}(\e) 
\prod_{j=1}^{N_{\e}(\omega)} p_{n+1-t^j_{\e}(\omega)}(e).
\end{aligned}
\ee

Given a path going from $\omega_0$ to $\omega_n$, all the edges $e$ in between $\omega_0$ and 
$\omega_n$ pointing in the direction of $\omega_n$, which we denote by $\cE(\omega_0,\omega_n)$, 
are traversed one time more than their reverse edges, while all other edges are traversed as 
often as their reverse edges. Therefore we obtain, assuming w.l.o.g.\ that $\omega_n>\omega_0$
(or $\theta>0$),
\be{Pshiftrel}
\begin{aligned}
\log \frac{P^{\xi}(\Theta\omega)}{P^{\xi}(\omega)} 
&= \sum_{e\in {\cE (\omega_0,\omega_n)}}
\log \frac{p_{n+1-t^{N_e(\omega)}_e(\omega)}(\e)}{p_{t^{N_e(\omega)}_e(\omega)}(e)}\\
&\qquad + \sum_{e\in E^+(\omega)}\sum_{j=1}^{N_{\e} (\omega)}
\log\left(\frac{p_{n+1-t^j_{\e}(\omega)}(e) p_{n+1-t^j_e(\omega)}(\e)}{p_{t^j_e(\omega)}(e)
p_{t^j_{\e}(\omega)}(\e)}\right).
\end{aligned}
\ee
In the \emph{static} random environment we have $p_n(e)=p_0(e)$ for all $n\in\N$ and $e\in 
E(\omega)$, and hence the second sum in (\ref{Pshiftrel}) is identically zero, while by the 
ergodic theorem the first sum equals
\be{ergstat}
(\omega_n-\omega_0) \langle\log[p_0(1,0)/p_0(0,1)]\rangle_{\nu_\rho} +o(n) 
= (\omega_n-\omega_0)(2\rho-1) \log (p/q) + o(n),
\qquad n\to\infty, 
\ee
where $\nu_\rho$ is the Bernoulli product measure on $\Omega$ with density $\rho$ (which is the 
law that is typically chosen for the static random environment). In the \emph{dynamic} random
environment, both sums in (\ref{Pshiftrel}) still look like ergodic sums, but since in general
\be{pnotid}
p_{t^j_e(\omega)}(e)\neq p_{t^i_e(\omega)}(e), \qquad i\neq j,
\ee
we have to use space-time ergodicity. 

\medskip\noindent
{\bf 2.\ Space-time ergodicity.}
Rewrite (\ref{Pshiftrel}) as
\be{Pshiftrel2}
\begin{aligned}
\log \frac{P^{\xi}(\Theta\omega)}{P^{\xi}(\omega)} 
&= \sum_{e\in {\cE (\omega_0,\omega_n)}} \log p_{n+1-t^{N_e(\omega)}_e(\omega)}(\e)
-\sum_{e\in {\cE(\omega_0,\omega_n)}} \log{p_{t^{N_e(\omega)}_e(\omega)}(e)}\\
&\qquad + \sum_{e\in E^+(\omega)} \log p_{n+1-t^1_{\e}(\omega)}(e)
+ \sum_{e\in E^+(\omega)} \log p_{n+1-t^1_e (\omega)}(\e)\\
&\qquad - \sum_{e\in E^+(\omega)} \log p_{t^1_e(\omega)}(e)
- \sum_{e\in E^+(\omega)} \log p_{t^1_{\e}(\omega)}(\e)\\
&\qquad + \sum_{e\in E^+(\omega)} \sum_{j=2}^{N_{\e}(\omega)}
\log\left(\frac{p_{n+1-t^j_{\e}(\omega)}(e)p_{n+1-t^j_e(\omega)}(\e)}{p_{t^j_e(\omega)}(e)
p_{t^j_{\e}(\omega)}(\e)}\right),
\end{aligned}
\ee
and note that all the sums in (\ref{Pshiftrel2}) are of the form
\be{sumform1} 
\sum^N_{i=1} \log p_{t(i)}(\omega_0+i)
= \left\{\begin{array}{ll} (\log p)\displaystyle\sum_{i=1}^N \xi_{t_i}(\omega_0+i)
+ (\log q)\displaystyle\sum_{i=1}^N[1-\xi_{t_i}(\omega_0+i)],\\
(\log q)\displaystyle\sum_{i=1}^N \xi_{t_i}(\omega_0+i)
+ (\log p)\displaystyle\sum_{i=1}^N[1-\xi_{t_i}(\omega_0+i)],
\end{array} 
\right.
\ee
where  $t_i=t((i,i+1))$, with $t=t(\omega)\colon\,\{0,1,\dots,N\} \to \{0,1,\dots,n\}$ either 
strictly increasing or strictly decreasing with image set $I_n(t) \subset \{0,1,\dots,n\}$ 
such that $|I_n(t)|$ is of order $n$. Note that $N=N(\omega)=|\cE(\omega_0,\omega_n)|=\omega_n
-\omega_0=\lfloor\theta n\rfloor$ in the first two sums in (\ref{Pshiftrel2}), $N=N(\omega)
=|E^+(\omega)|\geq\omega_n-\omega_0=\lfloor\theta n\rfloor$ in the remaining sums, and
\be{TimeBound}
|t_j-t_i| \geq j-i, \qquad j>i.
\ee

The aim is to show that 
\be{ergsum} 
\lim_{N\rightarrow\infty} \frac{1}{N} \sum^N_{i=1} 
\log p_{t_i}(i) = \langle \log p_{0}(0) \rangle_\mu 
= \rho \log p + (1-\rho) \log q \qquad \xi-a.s. \mbox{ for all } \omega
\ee
or, equivalently,
\be{ergsum2}
\lim_{N\rightarrow\infty} \frac{1}{N} \sum^N_{i=1} \xi_{t_i}(i) 
= \langle \xi_{0}(0) \rangle_\mu = \rho \qquad \xi-a.s. \mbox{ for all } \omega,
\ee 
where, since we take the limit $N\to\infty$, we think of $\omega$ as an infinite path in which
the $n$-step path $(\omega_0,\dots,\omega_n)$ with $\omega_0=0$ and $\omega_n=\lfloor \theta n
\rfloor$ is embedded. Because $P^\mu$ is tail trivial (recall (\ref{tailtriv})) and 
$\lim_{i\to\infty} t_i=\infty$ for all $\omega$ by (\ref{TimeBound}), the limit exists 
$\xi$-a.s.\ for all $\omega$. To prove that the limit equals $\rho$ we argue as follows. Write
\be{ergvar2} 
\mathrm{Var}^{P^\mu}\Bigg(\frac{1}{N} \sum^N_{i=1} \xi_{t_i}(i)\bigg)
=\frac{\rho(1-\rho)}{N}
+ \frac{2}{N^2} \sum_{i=1}^N \sum_{j>i}\mathrm{Cov}^{P^\mu}\big(\xi_{t_i}(i),\xi_{t_j}(j)\big).
\ee
By (\ref{erg}), we have
\be{corshift}
\mathrm{Cov}^{P^\mu}\big(\xi_{t_i}(i),\xi_{t_j}(j)\big)
= \mathrm{Cov}^{P^\mu}\big(\xi_0(0),\xi_{|t_j-t_i|}(j-i)\big).
\ee
In view of (\ref{TimeBound}), it therefore follows that
\be{corlim}
\lim_{k\to\infty} \sup_{l \geq k} \mathrm{Cov}^{P^\mu}\big(\xi_0(0),\xi_l(k)\big) = 0
\quad \Longrightarrow \quad 
\lim_{N\to\infty} \mathrm{Var}^{P^\mu}\Bigg(\frac{1}{N} \sum^N_{i=1} \xi_{t_i}(i)\bigg) = 0.
\ee
But the l.h.s.\ of (\ref{corlim}) is true by the tail triviality of $P^\mu$.

\medskip\noindent
{\bf 3.\ Implication for the rate function.}
Having proved (\ref{ergsum}) holds, we can now use (\ref{Pshiftrel2}--\ref{sumform1}) 
and (\ref{ergsum}--\ref{ergsum2}) to obtain
\be{ratio}
\frac{P^{\xi}(\Theta\omega)}{P^{\xi}(\omega)}
= \erom^{A(\omega_n-\omega_0)+o(n)} \quad \mbox{ with } \quad A=(2\rho-1)\log\,(p/q).
\ee 
Thus, the probability that the walk moves from $0$ to $\lfloor\theta n\rfloor$ in 
$n$ steps is given by
\be{Psumrel}
\begin{aligned}
&P^\xi(\omega_n=\lfloor\theta n\rfloor \mid \omega_0=0) 
= \sum_{ {\omega\colon\,|\omega|=n} \atop {\omega_0=0,\omega_n=\lfloor\theta n\rfloor} } 
P^{\xi}(\omega)
= \sum_{ {\omega\colon\,|\omega|=n} \atop {\omega_0=0,\omega_n=\lfloor\theta n\rfloor} } 
P^{\xi}(\Theta\omega)\,\erom^{-A\lfloor\theta n\rfloor+o(n)}\\
&= \erom^{-A\lfloor\theta n\rfloor+o(n)} 
\sum_{ {\omega\colon\,|\omega|=n} \atop {\omega_0=\lfloor\theta n\rfloor,\omega_n=0} } 
P^{\xi}(\omega) 
\qquad = \erom^{-A\lfloor\theta n\rfloor+o(n)}\,
P^{\xi}(\omega_n=0 \mid \omega_0=\lfloor\theta n\rfloor).
\end{aligned}
\ee 
Since the quenched rate function is $\xi$-a.s.\ constant, we have
\be{PIrel1}
\begin{aligned}
P^\xi(\omega_n=\lfloor\theta n\rfloor \mid \omega_0=0)
&= P_{0}^{\xi}(X_n=\lfloor\theta n\rfloor) = \erom^{-nI^\mathrm{que}(\theta)+o(n)},\\ 
P^{\xi}(\omega_n=0 \mid \omega_0=\lfloor\theta n\rfloor) 
&= P^{\tau_{\lfloor\theta n\rfloor}\xi}_{0}(X_n=-\lfloor\theta n\rfloor) 
= \erom^{-nI^\mathrm{que}(-\theta)+o(n)},
\end{aligned}
\ee 
and hence
\be{PIsym} 
\frac{1}{n} \log \left(\frac{P^{\xi}(\omega_n=\lfloor\theta n\rfloor \mid \omega_0=0)}
{P^{\xi}(\omega_n=0 \mid \omega_0=\lfloor\theta n\rfloor)} \right)
\to -I^\mathrm{que}(\theta)+I^\mathrm{que}(-\theta). 
\ee 
Together with (\ref{Psumrel}), this leads to the symmetry relation
\be{IquesymA}
-I^\mathrm{que}(\theta) + I^\mathrm{que}(-\theta) = -A\theta.
\ee

\medskip\noindent
{\bf 4.\ From discrete to continuous time.}
Let $\chi=(\chi_n)_{n\in\N_0}$ denote the jump times of the continuous-time random 
walk $X=(X_t)_{t \geq 0}$ (with $\chi_0=0$). Let $Q$ denote the law of $\chi$. The 
increments of $\chi$ are i.i.d.\ random variables, independent of $\xi$, whose 
distribution is exponential with mean $1/(\alpha+\beta)$. Define
\be{tildexiX}
\begin{array}{lllll}
\xi^* &=& (\xi^*_n)_{n\in\N_0} \quad \mbox{ with } \quad \xi^*_n &=& \xi_{\chi_n},\\ 
X^* &=& (X^*_n)_{n\in\N_0} \quad \mbox{ with} \quad X^*_n &=& X_{\chi_n}. 
\end{array}
\ee
Then $X^*$ is a discrete-time random walk in a random environment $\xi^*$ of the type 
considered in Steps 1--3, with $p=\alpha/(\alpha+\beta)$ and $q=\beta/(\alpha+\beta)$.
The analogue of (\ref{ergsum2}) reads
\be{ergsum2alt}
\lim_{N\to\infty} \frac{1}{N} \sum_{i=1}^N \xi_{\chi_{t_i}}(i) = \rho 
\qquad \xi,\chi-a.s. \mbox{ for all } \omega,
\ee
where we use that the law of $\chi$ is invariant under permutations of its increments.
All we have to do is to show that
\be{varalt}
\lim_{N\to\infty}
E^Q\left(\mathrm{Var}^{P^\mu}\left(\frac{1}{N} 
\sum_{i=1}^N \xi_{\chi_{t_i}}(i)\right)\right) = 0.
\ee
But 
\be{corlimprop}
E^Q\left(\mathrm{Cov}^{P^\mu}\left(\xi_{\chi_{t_i}}(i),\xi_{\chi_{t_j}}(j)\right)\right)
= E^Q\left(\mathrm{Cov}^{P^\mu}\left(\xi_0(0),\xi_{|\chi_{t_j}-\chi_{t_i}|}(j-i)\right)\right),
\ee
while (\ref{TimeBound}) ensures that $\lim_{j\to\infty} |\chi_{t_j}-\chi_{t_i}| \to \infty$ 
$\chi$-a.s.\ for all $\omega$ as $j-i\to\infty$. Together with the tail triviality of $P^\mu$
assumed in (\ref{tailtriv}), this proves (\ref{varalt}).   
\end{proof}


\section{Proof of Theorem~\ref{SSEslow}}
\label{S4}

In Section~\ref{S4.1} we show that the simple symmetric exclusion process suffers
traffic jams. In Section~\ref{S4.2} we prove that these traffic jams cause the
slow-down of the random walk. 


\subsection{Traffic jams}
\label{S4.1}

In this section we derive two lemmas stating that long strings of occupied and 
vacant sites have an appreciable probability to survive for a long time under 
the simple symmetric exclusion dynamics, both when they are alone (Lemma \ref{sselds}) 
and when they are together but sufficiently separated from each other (Lemma \ref{sseldsalt}). 
These lemmas, which are proved with the help of the \emph{graphical representation}, are 
in the spirit of \cite{Ar85}.

\begin{figure}[hbtp]
\vspace{1cm}
\begin{center}
\setlength{\unitlength}{0.3cm}
\begin{picture}(20,10)(0,0)
\put(0,0){\line(22,0){22}}
\put(0,11){\line(22,0){22}}
\put(2,0){\line(0,12){12}}
\put(5,0){\line(0,12){12}}
\put(8,0){\line(0,12){12}}
\put(11,0){\line(0,12){12}}
\put(14,0){\line(0,12){12}}
\put(17,0){\line(0,12){12}}
\put(20,0){\line(0,12){12}}
\qbezier[15](2.1,2)(3.5,2)(4.9,2)
\qbezier[15](5.1,4)(6.5,4)(7.9,4)
\qbezier[15](8.1,7)(9.5,7)(10.9,7)
\qbezier[15](2.1,8)(3.5,8)(4.9,8)
\qbezier[15](11.1,2.5)(12.5,2.5)(13.9,2.5)
\qbezier[15](17.1,4)(18.5,4)(19.9,4)
\qbezier[15](11.1,4.5)(12.5,4.5)(13.9,4.5)
\qbezier[15](14.1,7.5)(15.5,7.5)(16.9,7.5)
\put(10.7,-1.2){$x$}
\put(8.4,11.6){$y$}
\put(-1.2,-.3){$0$}
\put(-1.2,10.7){$t$}
\put(12,2.7){$\rightarrow$}
\put(12,4.7){$\leftarrow$}
\put(9,7.3){$\leftarrow$}
\put(11.3,.8){$\uparrow$}
\put(14.3,3){$\uparrow$}
\put(11.3,5.5){$\uparrow$}
\put(8.3,8.8){$\uparrow$}
\put(11,0){\circle*{.35}}
\put(8,11){\circle*{.35}}
\put(23,0){$\Z^d$}
\end{picture}
\end{center}
\caption{\small Graphical representation. The dashed lines are
links. The arrows represent a path from $(x,0)$ to $(y,t)$.}
\label{fig-graph}
\end{figure}
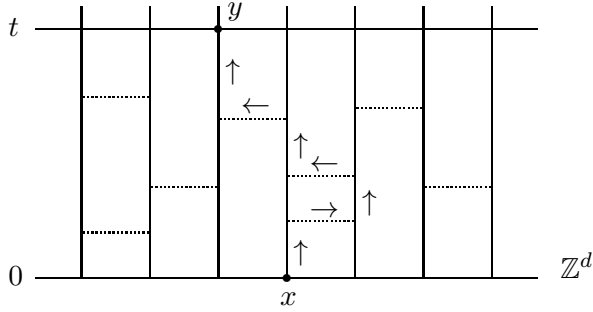

In the \emph{graphical representation} of the simple symmetric exclusion process, space is 
drawn sidewards, time is drawn upwards, and for each pair of nearest-neighbor sites $x,y\in\Z$ 
links are drawn between $x$ and $y$ at Poisson rate $1$. The configuration at time $t$ is 
obtained from the one at time $0$ by transporting the local states along paths that move 
upwards with time and sidewards along links (see Fig.~\ref{fig-graph}).


\bl{sselds}
There exists a $C=C(\rho)>0$ such that, for all $Q\subset\Z$ and all $t \geq 1$,
\be{sseldprop}
P^{\nu_\rho}\Big(\xi_s(x)=0\,\,\forall\,x\in Q\,\,\forall\,s\in [0,t]\Big) 
\geq \erom^{-C|Q|\sqrt{t}}.
\ee
\el

\begin{proof}
Let
\be{Hdef} 
H^Q_t = \Big\{x\in\Z\colon\,\,\exists \mbox{ path in }
\cG \mbox{ from } (x,0) \mbox{ to } Q \times [0,t]\Big\}.
\ee
Note that $H^Q_0=Q$ and that $t \mapsto H^Q_t$ is non-decreasing. Denote by $\cP$
and $\cE$, respectively, probability and expectation w.r.t.\ $\cG$. Let
$V_0=\{x\in\Z\colon\,\xi_0(x)=0\}$ be the set of initial
locations of the vacancies. Then
\be{Pssefull}
P^{\nu_\rho}\Big(\xi_s(x)=0\,\,\forall\, x\in Q\,\,\forall\, s\in [0,t]\Big) 
= (\cP\otimes\nu_\rho)\left(H^Q_t \subset V_0\right).
\ee
Indeed, if $\xi_0(x)=1$ for some $x\in H^Q_t$, then this $1$ will propagate into $Q$
prior to time $t$ (see Fig.~\ref{fig-graphappl}).


\begin{figure}[hbtp]
\vspace{1cm}
\begin{center}
\setlength{\unitlength}{0.3cm}
\begin{picture}(20,10)(0,0)
\put(0,0){\line(22,0){22}}
\put(0,11){\line(22,0){22}}
\put(2,0){\line(0,12){12}}
\put(5,0){\line(0,12){12}}
\put(8,0){\line(0,12){12}}
\put(11,0){\line(0,12){12}}
\put(14,0){\line(0,12){12}}
\put(17,0){\line(0,12){12}}
\put(20,0){\line(0,12){12}}
\qbezier[15](2.1,2)(3.5,2)(4.9,2)
\qbezier[15](5.1,4)(6.5,4)(7.9,4)
\qbezier[15](8.1,7)(9.5,7)(10.9,7)
\qbezier[15](2.1,8)(3.5,8)(4.9,8)
\qbezier[15](11.1,2.5)(12.5,2.5)(13.9,2.5)
\qbezier[15](17.1,4)(18.5,4)(19.9,4)
\qbezier[15](11.1,4.5)(12.5,4.5)(13.9,4.5)
\qbezier[15](14.1,7.5)(15.5,7.5)(16.9,7.5)
\put(1.7,-1.2){$x$}
\put(-1,-.3){$0$}
\put(-1,10.7){$t$}
\put(10,-1.5){$[$}
\put(20.9,-1.5){$]$}
\put(15.2,-1.5){$Q$}
\put(16.8,-1.5){$\longrightarrow$}
\put(12.5,-1.5){$\longleftarrow$}
\put(3,2.5){$\rightarrow$}
\put(6,4.5){$\rightarrow$}
\put(9,7.5){$\rightarrow$}
\put(2.5,.6){$\uparrow$}
\put(5.5,2.6){$\uparrow$}
\put(8.5,5.1){$\uparrow$}
\put(2,0){\circle*{.35}}
\put(11,7){\circle*{.35}}
\end{picture}
\end{center}
\caption{\small A path from $(x,0)$ to $Q \times [0,t]$.}
\label{fig-graphappl}
\end{figure}


By Jensen's inequality,
\be{sseJen}
(\cP\otimes\nu_\rho)\left(H^Q_t\subset V_0\right)
= \cE\left((1-\rho)^{|H^Q_t|}\right) \geq (1-\rho)^{\cE(|H^Q_t|)}.
\ee
Moreover, since $H^Q_t=\cup_{x\in Q} H^x_t$ and $\cE(|H^x_t|)$ does
not depend on $x$, we have
\be{HQtmean}
\cE(|H^Q_t|) \leq |Q|\,\cE(|H^0_t|),
\ee
and, by time reversal, we see that
\be{H0tmean}
\begin{aligned}
\cE(|H^0_t|) &= \sum_{x\in\Z} \cP\Big(\exists \mbox{ path in } \cG
\mbox{ from } (x,0) \mbox{ to } \{0\}\times [0,t]\Big)\\
&= \sum_{x\in\Z} \P_0^{\SRW}(\tau_x\leq t) =\E_0^{\SRW}(|R_t|),
\end{aligned}
\ee
where $\P_0^{\SRW}$ is the law of simple symmetric random walk jumping at rate $1$ 
starting from $0$, $R_t$ is the range (= number of distinct sites visited) at time 
$t$ and $\tau_x$ is the first hitting time of $x$. Combining (\ref{Pssefull}--\ref{H0tmean}),
we get
\be{sseldpropest}
P^{\nu_\rho}\Big(\xi_s(x)=0\,\,\forall\,x\in Q\,\,\forall\,s\in
[0,t]\Big) \geq (1-\rho)^{|Q|\,\E_0^{\SRW}(|R_t|)}.
\ee
The claim now follows from the fact that $R_0=1$ and $\E_0^{\SRW}(|R_t|) \sim C'\sqrt{t}$ 
as $t\to\infty$ for some $C'>0$ (see \cite{Sp76}, Section 1).
\end{proof}

\bl{sseldsalt} 
There exist $C=C(\rho)>0$ and $\delta>0$ such that, for all intervals $Q,Q'\subset\Z$ 
separated by a distance at least $2\sqrt{t\log t}$ and all $t\geq 1$, 
\be{sseldspropalt}
P^{\nu_\rho}\Big\{\xi_s(x)=1,\,\xi_s(y)=0\,\,\forall\,x\in
Q\,\,\forall\,y\in Q'\,\, \forall\,s\in [0,t]\Big\} 
\geq \delta\,\erom^{-C(|Q|+|Q'|)\sqrt{t}}. 
\ee 
\el

\begin{proof}
Recall (\ref{Hdef}) and abbreviate $A_t=\{H^Q_t \cap H^{Q'}_t = \emptyset\}$. Similarly as in 
(\ref{Pssefull}--\ref{sseJen}), we have 
\be{PAtrew} 
\mbox{l.h.s.}(\ref{sseldspropalt}) 
= (\cP \otimes \nu_{\rho})(A_t) 
= \cE\left(1_{A_t}\,\rho^{|H^Q_t|}\,(1-\rho)^{|H^{Q'}_t|}\right).
\ee 
Both $|H^{Q}_t|$ and $|H^{Q'}_t|$ are non-decreasing in the number of arrows in $\cG$, while 
$1_{A_t}$ is non-increasing in the number of arrows in $\cG$. Therefore, by the FKG-inequality 
(\cite{Li85}, Chapter II), we have 
\be{Aprodest}
\cE\left(1_{A_t}\,\rho^{|H^Q_t|}(1-\rho)^{|H^{Q'}_t|}\right) 
\geq \cP(A_t)\,\cE\Big(\rho^{|H^Q_t|}\Big)\,\cE\Big((1-\rho)^{|H^{Q'}_t|}\Big).
\ee 
We saw in the proof of Lemma \ref{sselds} that, for $t \geq 1$ and some $C>0$, 
\be{Htprods}
\cE\Big(\rho^{|H^Q_t|}\Big)\,\cE\Big((1-\rho)^{|H^{Q'}_t|}\Big)
\geq \erom^{-C(|Q|+|Q'|)\sqrt{t}}. 
\ee 
Thus, to complete the proof it suffices to show that there exists a $\delta>0$ 
such that 
\be{HxHy1} 
\cP(A_t) \geq \delta \mbox{ for } t \geq 1. 
\ee 

To that end, let $q=\max\{x\in Q\}$, $q'=\min\{x'\in Q'\}$ (where without loss of
generality we assume that $Q$ lies to the left of $Q'$). Then, using that $Q,Q'$ 
are intervals, we may estimate (see Fig.~\ref{fig-graphappl})
\be{HxHy3}
\begin{aligned}
\cP([A_t]^c) &= \cP\Big(\exists\,z\in \Z\colon\,(z,0) \to \partial
Q \times [0,t],\,
(z,0)\to \partial Q' \times [0,t]\Big)\\
&\leq \sum_{ {x\in\partial Q} \atop {x'\in\partial Q'} } \int_0^t
\Big[\cP\Big(\exists\,z\in \Z\colon\,(z,0) \to x\times [s,s+\di s],\,
(x,s)\to x'\times[s,t]\Big)\\
&\qquad\qquad\qquad + \cP\Big(\exists\,z\in \Z\colon\,(z,0) \to x'
\times [s,s+ds],\,
(x',s)\to x \times[s,t]\Big)\Big]\\
&= \sum_{ {x\in\partial Q} \atop {x'\in\partial Q'} } \int_0^t
\Big[\cP\Big(\exists\,z\in \Z\colon\,(z,0) \to x \times
[s,s+\di s]\Big)\,
\cP\Big((x,s)\to x' \times[s,t]\Big)\\
&\qquad\qquad\qquad + \cP\Big(\exists\,z\in \Z\colon\,(z,0) \to x'
\times [s,s+\di s]\Big)\,
\cP\Big((x',s)\to x \times[s,t]\Big)\Big]\\
&\leq 4\int_0^t \cP\Big(\exists\,z\in\Z\colon\, (z,0) \to 0 \times
[s,s+\di s]\Big)\,
\cP\Big((0,0) \to q'-q \times[0,t-s]\Big)\\[0.2cm]
&\leq 4\,\E_0^{\SRW}(|R_t|)\,\P_0^{\SRW}(\tau_{q'-q}\leq t),
\end{aligned}
\ee 
where the last inequality uses (\ref{H0tmean}). We already saw that $\,\E_0^{\SRW}(|R_t|) 
\sim C'\sqrt{t}$ as $t\to\infty$. By using, respectively, the reflection principle, 
the fact that $q'-q \geq 2\sqrt{t\log t}$, and the  Azuma-Hoeffding inequality (see 
\cite{Wi90}, (E14.2)), we get 
\be{srwbmrel} 
\P_0^{\SRW}(\tau_{q'-q}\leq t) = 2 \P_0^{\SRW}(S_t\geq {q'-q})
\leq 2 \P_0^{\SRW}(S_t\geq 2\sqrt{t\log t})
\leq 2 e^{-\frac{4t\log t}{2t}}=\frac{2}{t^2}.
\ee 
Combining (\ref{HxHy3}--\ref{srwbmrel}), we get $\cP([A_t]^c) \leq 2C'/t^{3/2}$, which 
tends to zero as $t\to\infty$. This proves the claim in (\ref{HxHy1}), because $\cP(A_t)>0$ 
for all $t \geq 0$.
\end{proof}


\subsection{Slow-down}
\label{S4.2}

We are now ready to prove Theorem~\ref{SSEslow}. The proof comes in two lemmas.

\bl{trapone}
For all $\rho\in (0,1)$ and $C>1/\log(\alpha/\beta)$,
\be{Flat}
\begin{aligned}
\lim_{t\to\infty} \frac{1}{t} \log\P_{\nu_\rho,0}(X_t\leq C\log t) &= 0,\\
\lim_{t\to\infty} \frac{1}{t} \log\P_{\nu_\rho,0}(X_t \geq -C\log t) &= 0.
\end{aligned}
\ee
\el

\begin{proof}
To prove the first half of (\ref{Flat}), the idea is to force $\xi$ to vacate an interval of 
length $C\log t$ to the right of $0$ up to time $t$ and to show that, with probability tending 
to $1$ as $t\to\infty$, $X$ does not manage to cross this interval up to time $t$ when $C$ is 
large enough.

For $t>0$, let $L_t=C\log t$ and
\be{Etdef} 
E_t = \big\{\xi_s(x)=0\,\,\forall\,x \in [0,L_t] \cap \Z\,\,\forall\,s\in [0,t]\big\}.
\ee
By Lemma \ref{sselds} we have, for some $C'>0$ and $t$ large enough,
\be{PEtest}
P^{\nu_\rho}(E_t) \geq \erom^{-C'\sqrt{t}\log t}. 
\ee
Hence
\be{PXtest} 
\P_{\nu_\rho,0}(X_t \leq L_t) 
\geq \P_{\nu_\rho,0}(X_t \leq L_t \mid E_t)\,P^{\nu_\rho}(E_t) 
\geq \P_{\nu_\rho,0}(X_t \leq L_t \mid E_t)\,\erom^{-C'\sqrt{t}\log t}. 
\ee
To complete the proof it therefore suffices to show that
\be{XtLtineq} 
\lim_{t\to\infty} \P_{\nu_\rho,0}(X_t \leq L_t \mid E_t) = 1. 
\ee

Let $\tau_{L_t}=\inf\{t\geq 0\colon\,X_t>L_t\}$. Then $\{X_t \leq L_t\mid E_t\} \supset 
\{\tau_{L_t}>t \mid E_t\}$, and so it suffices to show that
\be{XtLt3} 
\lim_{t\to\infty} \P_{\nu_\rho,0}(\tau_{L_t}>t \mid E_t) = 1. 
\ee
We say that $X$ starts a \emph{trial} when it enters the interval $[0,L_t] \cap \Z$ from 
the left prior. We say that the trial is successful when $X$ hits $L_t$ before returning 
to $0$. Let $M(t)$ be the number of trials prior to time $t$, and let $A_n$ be the event 
that the $n$-th trial is successful. Since
\be{tauAndef}
\{\tau_{L_t} \leq t\} \subset\bigcup_{n=1}^{M(t)} A_n, 
\ee
we have
\be{XtLt4}
\begin{aligned}
\P_{\nu_\rho,0}\big(\tau_{L_t} \leq t \mid E_t\big)
&\leq \P_{\nu_\rho,0}\left(\bigcup_{n=1}^{M(t)} A_n \,\Big|\, E_t\right)\\
&\leq \P_{\nu_\rho,0}\left(\bigcup_{n=1}^{2(\alpha+\beta)t} A_n,
\,M(t) \leq 2(\alpha+\beta)t \,\Big|\, E_t\right)\\
&\qquad + \P_{\nu_\rho,0}\Big(M(t)>2(\alpha+\beta)t \mid E_t\Big).
\end{aligned}
\ee
We will show that both terms in the r.h.s.\ tend to zero as $t\to\infty$. 

To estimate the second term in (\ref{XtLt4}), let $N(t)$ be the number of jumps by $X$ prior 
to time $t$, which is Poisson distributed with mean $(\alpha+\beta)t$ and is independent of 
$\xi$. Since $N(t)\geq M(t)$, it follows that
\be{Ntest}
\P_{\nu_\rho,0}\Big(M(t)>2(\alpha+\beta)t \mid E_t\Big) 
\leq \mathrm{Poi}\big(N(t)>2(\alpha+\beta)t\big), 
\ee
which tends to zero as $t\to\infty$. To estimate the first term in (\ref{XtLt4}), note
that, since $\P_{\nu_\rho,0}(A_n)$ is independent of $n$, we have
\be{Nest1}
\begin{aligned}
&\P_{\nu_\rho,0}\left(\bigcup_{n=1}^{2(\alpha+\beta)t} A_n,
\,M(t) \leq 2(\alpha+\beta)t \,\Big|\, E_t\right)\\
&\qquad \leq \P_{\nu_\rho,0}\left(\bigcup_{n=1}^{2(\alpha+\beta)t} A_n \,\Big|\, E_t\right) 
\leq 2(\alpha+\beta)t\,\,\P_{\nu_\rho,0}\big(A_1 \mid E_t\big).
\end{aligned}
\ee
But $\P_{\nu_\rho,0}\big(A_1 \mid E_t\big)$ is the probability that the random walk on $\Z$ 
that jumps to the right with probability $\beta/(\alpha+\beta)$ and to the left with
probability $\alpha/(\alpha+\beta)$ hits $L_t$ before $0$ when it starts from $1$.
Consequently,
\be{Nest2}
2(\alpha+\beta)t\,\,\P_{\nu_\rho,0}\big(A_1 \mid E_t\big)
= 2(\alpha+\beta)t\,\,\frac{(\alpha/\beta)-1}{(\alpha/\beta)^{L_t}-1},
\ee
which tends to zero as $t\to\infty$ when $L_t>C\log t$ with $C>1/\log(\alpha/\beta)$. 
This completes the proof of the first half of (\ref{Flat}).

To get the second half of (\ref{Flat}), note that $-X_t$ is equal in distribution to $X_t$ 
when $\rho$ is replaced by $1-\rho$.
\end{proof}

\bl{traptwo}
For all $\rho\in (0,1)$,
\be{Flattwo}
\lim_{t\to\infty} \frac{1}{t} \log\P_{\nu_\rho,0}(|X_t|\leq 2\sqrt{t\log t}) = 0.
\ee
\el

\begin{proof}
The idea is to create a trap around $0$ by forcing $\xi$ up to time $t$ to vacate an 
interval to the right of $0$ and occupy an interval to the left of $0$, separated by 
a suitable distance.


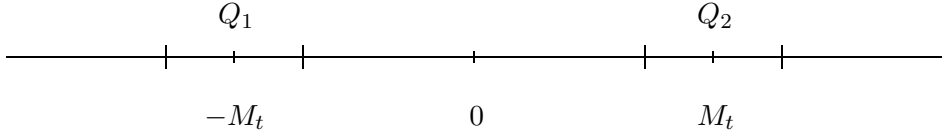
\begin{figure}[hbtp]
\vspace{-1.5cm}
\begin{center}
\setlength{\unitlength}{0.3cm}
\begin{picture}(12,10)(0,-2)
\put(-16,0){\line(41,0){41}} 
\put(-9,-0.5){\line(0,1){1}} 
\put(-6,-0.25){\line(0,1){0.5}}
\put(-3,-0.5){\line(0,1){1}} 
\put(4.5,-0.25){\line(0,1){0.5}}
\put(12,-0.5){\line(0,1){1}}
\put(15,-0.25){\line(0,1){0.5}}
\put(18,-0.5){\line(0,1){1}} 
\put(4.3,-3){$0$}
\put(-6.7,1.5){$Q_1$} 
\put(14.3,1.5){$Q_2$} 
\put(-7.3,-3){$-M_t$}
\put(14.3,-3){$M_t$}
\end{picture}
\end{center}
\caption{\small Location of the intervals $Q_1$ and $Q_2$. The width of $Q_1,Q_2$ is 
$2L_t$. The interval spanning $Q_1$, $Q_2$ and the space in between is $I_t$.} 
\label{fig-Qsets}
\end{figure}

For $t>0$, let $L_t=C\log t$ with $C>\log(\alpha/\beta)$, $M_t=\sqrt{t\log t}$,
\be{Qdefs}
Q_1 = \big(-M_t + [-L_t,L_t]\big) \cap\Z, \qquad 
Q_2 = \big(M_t + [-L_t,L_t]\big) \cap\Z,
\ee 
and $I_t = [-M_t- L_t,M_t + L_t] \cap\Z$ (see Fig.~\ref{fig-Qsets}). For $i=1,2$ and $j=0,1$, 
define the event 
\be{flat1} 
E_i^j = \Big\{\xi_s(x)=j\,\,\forall\,x\in Q_i,\,\,\forall\,s\in [0,t]\Big\}.
\ee 
Estimate, noting that $L_t \leq M_t$ for $t$ large enough, 
\be{flat2}
\begin{aligned}
&\P_{\nu_\rho,0}\Big(|X_t| \leq 2M_t\Big)
\geq \P_{\nu_\rho,0}\Big(X_t \in I_t\Big)\\
&\qquad \geq \P_{\nu_\rho,0}\Big(X_t \in I_t,\,E^1_1,E^0_2\Big)
=\P_{\nu_\rho,0}\Big(X_t \in I_t \mid E^1_1,E^0_2\Big)\,
\P_{\nu_\rho,0}\big(E^1_1,E^0_2\big).
\end{aligned}
\ee 
Since $\lim_{t\to\infty} \frac{1}{t} \log \P_{\nu_\rho,0}\big(E^1_1,E^0_2\big) = 0$
by Lemma~\ref{sseldsalt}, it suffices to show that 
\be{flat4}
\lim_{t\to\infty}
\P_{\nu_\rho,0}\Big(X_t \in I_t \mid E^1_1,E^0_2\Big) = 1.
\ee 
To that end, estimate
\be{1stcl}
\begin{aligned}
\P_{\nu_\rho,0}\Big(X_t \in I_t \mid E^1_1,E^0_2\Big)
&\geq \P_{\nu_\rho,0}\Big(X_t \leq M_t+L_t \mid E^1_1,E^0_2\Big)\\ 
&\qquad+ \P_{\nu_\rho,0}\Big(X_t \geq -M_t-L_t \mid E^1_1,E^0_2\Big) - 1.
\end{aligned}
\ee 
Now, irrespective of what $\xi$ does in between $Q_1$ and $Q_2$ up to time $t$, the same 
argument as in the proof of Lemma~\ref{trapone} shows that
\be{prlim12}
\begin{aligned}
&\lim_{t\to\infty} \P_{\nu_\rho,0}\Big(X_t \leq M_t + L_t \mid E^1_1,E^0_2\Big) = 1,\\
&\lim_{t\to\infty} \P_{\nu_\rho,0}\Big(X_t \geq -M_t-L_t \mid E^1_1,E^0_2\Big) = 1.
\end{aligned}
\ee
Combine this with (\ref{1stcl}) to obtain (\ref{flat4}). 
\end{proof}


\end{document}